\documentclass[aap,MSNbibl,nameyear,rotating,dvips]{arximspdf}
\usepackage{colortbl}
\usepackage{dcolumn}
\usepackage{graphicx}

\doi{10.1214/09-AAP653}
\volume{20}
\issue{3}
\pubyear{2010}
\firstpage{1098}
\lastpage{1125}

\makeatletter
\newcolumntype{d}[1]{D{.}{.}{#1}}
\newtheorem{theorem}{Theorem}
\newproclaim{conjecture}{Conjecture}
\newproclaim{remark}{Remark}
\makeatother

\begin{document}
\begin{frontmatter}

\title{A Markovian slot machine and Parrondo's paradox}
\runtitle{A Markovian slot machine}

\begin{aug}
\author[A]{\fnms{S. N.} \snm{Ethier}\ead[label=e1]{ethier@math.utah.edu}\corref{}} \and
\author[B]{\fnms{Jiyeon} \snm{Lee}\ead[label=e2]{leejy@yu.ac.kr}\thanksref{t1}}
\runauthor{S. N. Ethier and J. Lee}
\affiliation{University of Utah and Yeungnam University}
\address[A]{Department of Mathematics\\ University of Utah\\ 155 S.
1400 E.\\ Salt Lake City, Utah 84112\\ USA\\ \printead{e1}} 
\address[B]{Department of Statistics\\ Yeungnam University\\ 214-1
Daedong, Kyeongsan\\ Kyeongbuk 712-749\\ South Korea\\ \printead{e2}}
\thankstext{t1}{Supported by the Yeungnam University research grants
in 2008.}
\end{aug}

\received{\smonth{7} \syear{2009}}

%
\begin{abstract}
The antique Mills Futurity slot machine has two unusual features.
First, if a player loses 10 times in a row, the 10 lost coins are
returned. Second, the payout distribution varies from coup to coup in a
manner that is nonrandom and periodic with period 10. It follows that
the machine is driven by a 100-state irreducible period-10 Markov
chain. Here, we evaluate the stationary distribution of the Markov
chain, and this leads to a strong law of large numbers and a central
limit theorem for the sequence of payouts.
Following a suggestion of Pyke [In \textit{Mathematical Statistics and
Applications: Festschrift for Constance van Eeden} (2003)
185--216 Institute of Mathematical Statistics], we address the question of
whether there exists a two-armed version of this ``one-armed bandit''
that obeys Parrondo's paradox. More precisely, is there such a machine
with the property that the casino can honestly advertise that both arms
are fair, yet when players alternate arms in certain random or
nonrandom ways, the casino makes money in the long run? The answer is a
qualified yes. Although this ``history-dependent'' game is conceptually
simpler than the original such games of Parrondo, Harmer and Abbott
[\textit{Phys. Rev. Lett.} (2000) \textbf{85} 5226--5229], it is nearly
as complicated analytically, and open problems remain.
\end{abstract}

%

\end{frontmatter}
%

\section{Introduction}\label{sec1}

The Futurity slot machine, a 1936 design of Mills Novelty Company of
Chicago, has two unusual features, one readily apparent and the other
less so. The readily apparent feature is that, if the player loses 10
times in a row, the 10 lost coins are returned. At the top of the
machine is a pointer that indicates the number of consecutive losses
incurred. It advances by 1 after each loss, and resets at 0 after a win
or after 10 consecutive losses. The less apparent feature is that there
are 20 symbols on each of the three reels but only the ones in
even-numbered positions can appear on the payline if the machine is in
mode E, while only the ones in odd-numbered positions can appear on the
payline if the machine is in mode O. The mode is nonrandom and is
determined by a cam that rotates through 10 positions, advancing one
position with each coup and resulting in a specific mode pattern of
length 10, EEEEEOEEEO, which is repeated ad infinitum. (Note that we
could substitute any cyclic permutation of this mode pattern, such as
EEEOEEEEEO, without effect.) When in mode E, the machine is extremely
``tight'' (i.e., the mean payout from a one-coin bet is much less than
1). When in mode O, it is extremely ``loose.''

There are several questions that might be asked. Does the sequence of
payouts obey the strong law of large numbers and the central limit
theorem, as it would for a traditional slot machine for which the
sequence can be assumed independent and identically distributed? If so,
what are the mean and variance parameters? What is the asymptotic
probability of a nonzero payout? How frequently does the player lose 10
times in a row, thereby receiving the so-called Futurity award? Are
there advantageous opportunities depending on the information available
to the player about the state of the machine?

Notice that the machine is driven by a Markov chain with state space
$\Sigma:=\{0,1,\ldots,9\}\times\{0,1,\ldots,9\}$ interpreted as
follows. The machine is in state $(i,j)$ if the cam position is $i$ and
the pointer position is $j$. (If the cam position is 5 or 9, the
machine is in mode O; if the cam position is 0--4 or 6--8, the machine
is in mode~E.) If we kept track of the mode (E or O) instead of the cam
position (0--9), we would lose the Markov property. There is also a
pointer position 10, but from that position the pointer instantly moves
to position 0, so we can ignore pointer position 10. By evaluating the
stationary distribution of this Markov chain, we can infer the
long-term behavior of the slot machine. Specifically, we can establish
a strong law of large numbers and a central limit theorem for the
sequence of payouts.

The Futurity came to our attention via articles of Geddes (\citeyear{G}) and
Geddes and Saul (\citeyear{GS}) that appeared in \textit{Loose Change}, a
magazine for collectors of antique slot machines (published 1977--1998
and archived at the UNLV Lied Library). Geddes and Saul used Monte
Carlo simulation to study the Futurity, claiming that an analytical
solution ``falls somewhere between formidable and monumental on a
relative scale of mathematical difficulty.'' As we will see, the claim
is untrue.

Parrondo's paradox can be regarded as the observation that there exist
two fair games that can be combined, by either random mixture or
nonrandom alternation, to create an unfair game. See the survey
articles by Harmer and Abbott (\citeyear{HA}), Parrondo and Din\'is (\citeyear{PD}),
Epstein (\citeyear{E}) and Abbott (\citeyear{A}). To motivate his discussion of the
paradox, Pyke (\citeyear{P}) raised the following question without providing an
explicit answer.

\begin{quote}
You are about to play a two-armed slot machine. The casino that owns
this two-armed bandit advertises that both arms on their two-armed
machines are ``fair'' in the sense that any player who plays either of
the arms is assured that the average cost per play approaches zero as
the number of plays increases. However, the casino does not constrain
you to stay with one arm; you are allowed to use either arm on every
play. [ \dots] The question of interest in this context would be
whether it is possible for the casino to still make money using only
``fair'' games.
\end{quote}

Our aim here is to formulate a two-armed version of the Mills Futurity
that answers Pyke's question affirmatively. The feature of the Futurity
that permits Parrondian behavior is the Futurity award (the return of
the 10 lost coins after 10 consecutive losses); the periodicity of the
payout distribution is not important. This ``history-dependent'' bonus
feature makes our hypothetical two-armed slot machine not unlike the
history-dependent games introduced by Parrondo, Harmer and Abbott
(\citeyear{PHA}). In fact, it has some advantages over the original such games:
It is conceptually simpler and less contrived. On the other hand, it is
nearly as complicated analytically.

Actually, our answer to Pyke's question must be qualified. It is an
unqualified yes for the random-mixture strategies. It is a qualified
yes for the nonrandom-alternation strategies because certain
assumptions are needed and our conclusions rely on an unproved
conjecture. And the answer is simply no if the player's strategy is
completely unrestricted because there exist strategies that actually
give the player an advantage. In particular, our two-armed version of
the Futurity is not ready for casino play.

We should clarify how it works. The player can pull either arm at each
coup. After 10 consecutive losses, regardless of the order of play of
the two arms, the 10 lost coins are returned to the player. On the
other hand, each arm has its own cam mechanism, each with 10 positions,
hence its own periodic pattern of payout distributions (though the
payout distribution need not vary). The cam position for an arm
advances only when that arm is pulled. Indeed, if this were not the
case and both cam positions advanced when either arm was pulled, astute
players would simply pull the arm with the higher mean payout, and the
casino would be beaten at its own game. Of course, there is nothing
special about the number 10 in this context, so we replace it
throughout by the integer $J\ge2$.

The question of whether Parrondo's paradox can appear in the casino
setting was raised by Harmer and Abbott (\citeyear{HA}), Section 2.3.3. Our
example shows that the potential exists, even though it will not likely
be realized. However, in our case the winning game created from two
fair games is winning for the casino, not for the player. If it were
the other way around, the casino would likely discontinue the game or
change the rules.

In a previous paper [Ethier and Lee (\citeyear{EL})], the authors formulated a
general version of Parrondo's games. The results of that paper do not
immediately apply here because the present underlying irreducible
Markov chain is periodic. Even if that issue could be overcome, the
Markov chain here is rather complicated relative to the three- and
four-state chains that were studied in the previous paper. It is
therefore preferable to use a different approach here that avoids
having to evaluate the fundamental matrix and spectral representation
associated with the one-step transition matrix of the Markov chain.

\section{The Markov chain at equilibrium}\label{sec2}

We will analyze a generalized (one-armed) version of the Futurity,
dependent on several parameters. In Section \ref{sec5}, we will substitute the
actual numbers.

We assume that the cam controlling the payout distribution has $I$
positions, denoted by $0,1,\ldots,I-1$. When in cam position $i$, the
probability of a nonzero payout is $p_i$, the mean payout is $\mu_i$
and the variance of the payout is $\sigma_i^2$; none of these
parameters takes the Futurity award into account. As for the Futurity
award, we assume that, if the player loses $J$ times in a row, the $J$
lost coins are returned. A pointer that indicates the number of
consecutive losses advances by 1 after each loss, and resets at 0 after
a win or after $J$ consecutive losses.

If we were interested solely in the Futurity, we would take $I=J$ and
simplify matters considerably. However, in studying Parrondo's paradox
for a two-armed version of the Futurity, it will be necessary to allow
$I$ in the generalized one-armed machine to be an integer multiple of
$J$, say $I=dJ$ for a positive integer $d$. Of course, the case $d=1$
is included and is in fact of primary interest.

The Markov chain $\{(X_n,Y_n)\}_{n\ge0}$ that drives (or controls) the
generalized (one-armed) Futurity has state space $\Sigma:=\{0,1,\ldots
,I-1\}\times\{0,1,\ldots,J-1\}$. It is in state $(i,j)$ at time $n$
if the cam position is $i$ and the pointer position is $j$ following
the $n$th coup. The transition probabilities have a very simple form:
\begin{eqnarray*}
P((i,j),(k,l))&:=&\mathrm{P}\bigl((X_{n+1},Y_{n+1})=(k,l)\mid
(X_n,Y_n)=(i,j)\bigr)\\
& =&
\cases{
p_i&\quad if $(k,l)=\bigl(i+1\mbox{ } (\operatorname{mod} I),0\bigr)$ and $j\le
J-2$,\vspace*{2pt}\cr
q_i&\quad if $(k,l)=\bigl(i+1\mbox{ } (\operatorname{mod} I),j+1\bigr)$ and $j\le
J-2$,\vspace*{2pt}\cr
1&\quad if $(k,l)=\bigl(i+1\mbox{ } (\operatorname{mod} I),0\bigr)$ and $j=J-1$,}
\end{eqnarray*}
where $0<p_i<1$ and $q_i:=1-p_i$ for $i=0,1,\ldots,I-1$. We notice
that the one-step transition matrix $\mathbf{P}$ is irreducible and periodic
with period $I$.

\begin{theorem}\label{teo1}
The unique stationary distribution $\bolds{\pi}$ for the Markov chain in
$\Sigma$ with one-step transition matrix $\mathbf{P}$ is given recursively by
%
\begin{eqnarray}\label{1}
 \quad &&\pi(i,0)\nonumber\\[-8pt]\\[-8pt]
 \quad &&\qquad ={p_{i-1}+q_{i-1}\cdots q_{i-J}p_{i-J-1}+\cdots
+q_{i-1}\cdots q_{i-(d-1)J}p_{i-(d-1)J-1}\over I(1-Q)}\nonumber
\end{eqnarray}
for $i=0,1,\ldots,I-1$,
%
\begin{eqnarray}\label{2}
\pi(i,1)&=&q_{i-1}\pi(i-1,0),\qquad   i=0,1,\ldots,I-1,\\ \label{3}
\pi(i,2)&=&q_{i-1}\pi(i-1,1),\qquad   i=0,1,\ldots,I-1,\\
&\vdots& \nonumber\\ \label{4}
\pi(i,J-1)&=&q_{i-1}\pi(i-1,J-2),\qquad   i=0,1,\ldots,I-1,
\end{eqnarray}
where $Q:=q_0q_1\cdots q_{I-1}$, $p_{-i}:=p_{I-i}$ and
$q_{-i}:=q_{I-i}$ for $i=1,2,\ldots,I$, and $\pi(-1,j):=\pi(I-1,j)$
for $j=0,1,\ldots,J-1$. Furthermore,
%
\begin{equation}\label{5}
\pi(i,0)+\pi(i,1)+\cdots+\pi(i,{J-1})={1\over I},\qquad   i=0,1,\ldots,I-1.
\end{equation}
\end{theorem}

\begin{remark*}
In the special case $I=J$ (i.e., $d=1$), (\ref{1}) and (\ref{4})
simplify to
\[
\pi(i,0)={p_{i-1}\over J(1-Q)},\qquad  \pi(i,J-1)={p_i Q\over q_i J(1-Q)}.
\]
\end{remark*}

\begin{pf*}{Proof of Theorem \ref{teo1}}
The stationary distribution is the unique probability (row) vector $\bolds{\pi}$ satisfying
%
\begin{equation}\label{pi=piP}
\bolds{\pi}=\bolds{\pi}\mathbf{P}.
\end{equation}
Equations (\ref{2})--(\ref{4}) are immediate from this. This reduces
the problem to a system of $I$ linear equations in $I$ variables, $\pi
(i,0)$, $i=0,1,\ldots,I-1$. The system is a rather complicated one, so
we take a different approach, noticing that these probabilities can be
obtained probabilistically.

If the Markov chain has the stationary distribution as its initial
distribution, it is a stationary process, and we can extend its time
parameter to the set of all integers. Intuitively, we can assume that
the machine has been operating forever. What is the probability that,
at a particular time, the Markov chain is in state $(i,0)$? First the
cam position must be $i$, the probability of which is $1/I$. Second,
either the last coup resulted in a win (conditional probability
$p_{i-1}$) or the last coup completed a string of $J$ or $2J$ or $3J$
or \dots\ consecutive losses, causing the pointer to reset at 0 and
the Futurity award to be paid. Thus, the conditional probability that
the pointer position is 0, given that the cam position is $i$, is
\begin{eqnarray*}
&&p_{i-1}+q_{i-1}\cdots q_{i-J}p_{i-J-1}+q_{i-1}\cdots
q_{i-2J}p_{i-2J-1}+\cdots\nonumber\\
&&\quad   {}+q_{i-1}\cdots
q_{i-dJ}p_{i-dJ-1}+q_{i-1}\cdots q_{i-(d+1)J}p_{i-(d+1)J-1}+\cdots
\nonumber\\
&&\qquad  =p_{i-1}(1+Q+Q^2+\cdots)+q_{i-1}\cdots
q_{i-J}p_{i-J-1}(1+Q+Q^2+\cdots)+\cdots\nonumber\\
&& \quad \qquad  {}+q_{i-1}\cdots
q_{i-(d-1)J}p_{i-(d-1)J-1}(1+Q+Q^2+\cdots)\nonumber\\
&& \qquad {}=\bigl(p_{i-1}+q_{i-1}\cdots q_{i-J}p_{i-J-1}+\cdots\\
&&\quad \qquad\hspace*{10pt} {}+q_{i-1}\cdots q_{i-(d-1)J}p_{i-(d-1)J-1}\bigr)/(1-Q),
\end{eqnarray*}
where $p_{i-mI}:=p_i$ for all $i\in\{0,1,\ldots,I-1\}$ and $m\ge1$,
and similarly for complementary probabilities $q_{i-mI}$.
This implies (\ref{1}).

This argument is a bit heuristic [since we essentially assumed (\ref
{5}), one of the conclusions of the theorem], but now we can make it
rigorous. First,
we verify that $\bolds{\pi}$, given by (\ref{1})--(\ref{4}), is a
probability vector by proving (\ref{5}). Using (\ref{2})--(\ref{4})
and then~(\ref{1}),
the left-hand side of (\ref{5}) is equal to
\begin{eqnarray*}
&&\pi(i,0)+q_{i-1}\pi(i-1,0)+q_{i-1}q_{i-2}\pi(i-2,0)+\cdots\nonumber\\
&&\quad   {}+q_{i-1}\cdots q_{i-J+1}\pi(i-J+1,0)\nonumber
\\
&&\qquad  =\bigl[p_{i-1}+q_{i-1}\cdots q_{i-J}p_{i-J-1}+\cdots\\
&& \quad \qquad\hspace*{2pt}  {}+q_{i-1}\cdots q_{i-(d-1)J}p_{i-(d-1)J-1}\nonumber\\
&& \quad \qquad\hspace*{2pt}  {}+q_{i-1}\bigl(p_{i-2}+q_{i-2}\cdots
q_{i-J-1}p_{i-J-2}+\cdots\nonumber\\
&& \hspace*{76pt}   {}+q_{i-2}\cdots
q_{i-(d-1)J-1}p_{i-(d-1)J-2}\bigr)+\cdots\nonumber\\
&&\quad \qquad   {}+q_{i-1}\cdots q_{i-J+1}(p_{i-J}+q_{i-J}\cdots
q_{i-2J+1}p_{i-2J}+\cdots\nonumber\\
&&  \hspace*{134pt}\hspace*{-19pt}\hspace*{49pt}     {}
+q_{i-J}\cdots q_{i-dJ+1}p_{i-dJ})\bigr]/[I(1-Q)]\nonumber\\
&&\qquad  {}={p_{i-1}+q_{i-1}p_{i-2}+q_{i-1}q_{i-2}p_{i-3}+\cdots
+q_{i-1}\cdots q_{i-dJ+1}p_{i-dJ}\over I(1-Q)}\nonumber\\
&&\qquad  {}={1-q_{i-1}\cdots q_{i-dJ}\over I(1-Q)}={1\over I},
\end{eqnarray*}
where the second equality amounts to a rearrangement of terms, and the
third equality is an algebraic identity.

Next, for (\ref{pi=piP}) it will suffice to show, for $i=0,1,\ldots
,I-1$, that
\[
\pi(i,0)=p_{i-1}[\pi(i-1,0)+\cdots+\pi(i-1,J-2)]+\pi(i-1,J-1).
\]
This can be rewritten, using (\ref{5}) and (\ref{2})--(\ref{4}), as
\begin{eqnarray*}
\pi(i,0)&=&p_{i-1}[\pi(i-1,0)+\cdots+\pi(i-1,J-1)]+q_{i-1}\pi
(i-1,J-1)\nonumber\\
&=&{p_{i-1}\over I}+q_{i-1}\cdots q_{i-J}\pi(i-J,0).
\end{eqnarray*}
Fix $i$ and substitute (\ref{1}). It is enough that
\begin{eqnarray*}
&&p_{i-1}+q_{i-1}\cdots q_{i-J}p_{i-J-1}+\cdots+q_{i-1}\cdots
q_{i-(d-1)J}p_{i-(d-1)J-1}\nonumber\\
&&\qquad  =(1-Q)p_{i-1}+q_{i-1}\cdots q_{i-J}
(p_{i-J-1}+q_{i-J-1}\cdots q_{i-2J}p_{i-2J-1}+\cdots\\
&& \hspace*{50pt}\hspace*{171pt} {}+q_{i-J-1}\cdots q_{i-dJ}p_{i-dJ-1}).
\nonumber
\end{eqnarray*}
Canceling like terms, this reduces to $p_{i-1}=(1-Q)p_{i-1}+Qp_{i-1}$,
which proves that $\bolds{\pi}$, defined by (\ref{1})--(\ref{4}), is the
stationary distribution for $\mathbf{P}$.
\end{pf*}

At equilibrium, what is the probability $p^\circ$ that, at a
particular coup, the player wins the $J$-coin Futurity award by losing
for the $J$th (or $2J$th or $3J$th or \dots) consecutive time? This
happens if and only if the Markov chain is in state $(i,J-1)$ for some
$i\in\{0,1,\ldots,I-1\}$ just before the specified coup \textit{and}
that coup results in a loss. Using (\ref{1})--(\ref{4}), the
probability is
%
\begin{eqnarray}\label{p^circ}
p^\circ=\sum_{i=0}^{I-1} \pi(i,J-1)q_i={1\over I(1-Q)}\sum
_{i=0}^{I-1}\sum_{k=1}^d q_i\cdots q_{i-kJ+1}p_{i-kJ}.
\end{eqnarray}
Notice that the last of the $d$ terms in the inner sum is $Qp_i$.

Therefore the mean payout, at equilibrium, is
%
\begin{equation}\label{mu^*}
\mu^*:={1\over I}\sum_{i=0}^{I-1}\mu_i+Jp^\circ.
\end{equation}
Incidentally, in the special case $I=J$ (i.e., $d=1$), (\ref{p^circ})
reduces to
%
\begin{equation}\label{p^circ:d=1}
p^\circ=\Biggl({1\over J}\sum_{i=0}^{J-1}p_i\Biggr){Q\over1-Q}.
\end{equation}

\section{Strong law of large numbers}\label{sec3}

Mean payout is the most important statistic of a slot machine. It can
be interpreted as the long-term proportion of coins played that are
paid out to the player. The justification of this interpretation is the
strong law of large numbers, which is well known to hold for
traditional machines, whose sequence of payouts is independent and
identically distributed (i.i.d.). Does the same conclusion hold for the
Futurity, even though the independence assumption and the identically
distributed assumption fail?

We will show that the answer is affirmative.

Let $R_1,R_2,\ldots$ be the sequence of payouts of the slot machine
\textit{excluding} the Futurity awards, given that the initial state
$(X_0,Y_0)=(i_0,j_0)\in\Sigma$ is specified. This sequence clearly
satisfies the strong law of large numbers. Indeed, $R_1,R_2,\ldots$
are independent, uniformly bounded, nonnegative random variables, with
$\{R_{n+mI},  m\ge0\}$ identically distributed as the payout
distribution in cam position $i$ (which has mean $\mu_i$), where
$i_0+n-1\equiv i$ (mod $I$). We conclude that, if $n$ is a multiple of
$I$, then
\[
n^{-1}\mathrm{E}[R_1+\cdots+R_n]={1\over I}\sum_{i=0}^{I-1}\mu
_i=:\mu.
\]
It follows from a version of the strong law of large numbers for
independent, but not identically distributed, random variables that
\[
n^{-1}(R_1+\cdots+R_n)\to\mu\qquad  \mbox{a.s.}
\]

Now, how does this change when the Futurity awards are taken into account?
Let $R_1^*,R_2^*,\ldots$ be the sequence of payouts of the slot
machine \textit{including} the Futurity awards, given that the initial
state $(X_0,Y_0)=(i_0,j_0)\in\Sigma$ is specified. Notice that, for
each $n\ge1$, $Y_n$ is a nonrandom function of $(X_0,Y_0)$ and $1_{\{
R_1=0\}},\ldots,1_{\{R_n=0\}}$; in particular, $Y_{n-1}$ is
independent of $R_n$. Clearly,
\begin{eqnarray*}
R_n^*&=&R_n+J\cdot1_{\{Y_{n-1}=J-1, R_n=0\}}\nonumber\\
&=&R_n+J\sum_{i=0}^{I-1} 1_{\{(X_{n-1},Y_{n-1})=(i,J-1), R_n=0\}
}, \qquad  n\ge1.
\end{eqnarray*}
It follows that
\begin{eqnarray*}
{R_1^*+\cdots+R_n^*\over n}&=&{R_1+\cdots+R_n\over n}+J\sum
_{i=0}^{I-1}{1\over n}\sum_{l=1}^n 1_{\{(X_{l-1},Y_{l-1})=(i,J-1),
R_l=0\}}\nonumber\\
&\to&\mu+J\sum_{i=0}^{I-1}\pi(i,J-1)q_i=\mu+Jp^\circ=\mu^*
\qquad \mbox{a.s.},
\end{eqnarray*}
where $\mu^*$ is as in (\ref{mu^*}); here the limit assertion
requires additional justification. Since the Markov chain is finite,
irreducible and periodic,
\[
{1\over n}\sum_{l=1}^n1_{\{(X_{l-1},Y_{l-1})=(i,J-1)\}}\to\pi
(i,J-1)\qquad  \mbox{a.s.}
\]
for $i=0,1,\ldots,I-1$, hence
%
\begin{eqnarray}\label{ratio}
 \quad &&{1\over n}\sum_{l=1}^n1_{\{(X_{l-1},Y_{l-1})=(i,J-1), R_l=0\}
}\nonumber\\
&&\qquad  {}=\Biggl({1\over n}\sum_{l=1}^n1_{\{
(X_{l-1},Y_{l-1})=(i,J-1)\}}\Biggr)\biggl({\sum_{l=1}^n1_{\{
(X_{l-1},Y_{l-1})=(i,J-1), R_l=0\}}\over\sum_{l=1}^n1_{\{
(X_{l-1},Y_{l-1})=(i,J-1)\}}}\biggr)\\
&&\qquad  {}\to\pi(i,J-1) q_i\qquad  \mbox{a.s.} \nonumber
\end{eqnarray}
for $i=0,1,\ldots,J-1$. We are using the fact that the ratio of sums
in (\ref{ratio}) represents the proportion of visits to $(i,J-1)$
(through time $n-1$) that result in a Futurity award. At each visit to
$(i,J-1)$ the probability of such an award is $q_i$ and the results are
determined independently; hence the ratio tends to $q_i$ a.s.\ by the
strong law of large numbers.

We have established the following version of the strong law of large numbers.

\begin{theorem}\label{teo2}
Let $R_1^*,R_2^*,\ldots$ be the sequence of payouts of the generalized
Futurity slot machine starting in an arbitrary initial state
$(X_0,Y_0)=(i_0,j_0)$. Then
\[
n^{-1}(R_1^*+\cdots+R_n^*)\to\mu^*\qquad  \mbox{a.s.}
\]
\end{theorem}

Observe that we can similarly obtain the asymptotic frequency of
nonzero payouts, the so-called ``hit frequency'' (usually reported as a
percentage):
\begin{eqnarray}\label{p^*}
&&n^{-1}\bigl(1_{\{R_1^*>0\}}+\cdots+1_{\{R_n^*>0\}}\bigr)\nonumber\\
&&\qquad  {}={1_{\{R_1>0\}}+\cdots+1_{\{R_n>0\}}\over n}\nonumber\\[-8pt]\\[-8pt]
&& \qquad \quad  {}+{1_{\{Y_0=J-1, R_1=0\}}+\cdots+1_{\{Y_{n-1}=J-1,
R_n=0\}}\over n}\nonumber\\
&& \qquad {}\to{1\over I}\sum_{i=0}^{I-1}p_i+p^\circ=:p^*\qquad  \mbox{a.s.} \nonumber
\end{eqnarray}

\section{Central limit theorem}

The second-most important statistic of a slot machine is the variance
of the payout. (This is arguable. Some would say that the hit frequency
$p^*$ is more important.) The variance permits determination of the
asymptotic distribution of the cumulative number of coins paid out by
the machine, via the central limit theorem. The central limit theorem
is well known to hold for traditional machines, whose sequence of
payouts is i.i.d. Does the same conclusion hold for the Futurity, even
though the independence assumption and the identically distributed
assumption fail?

We will show in two steps that the answer is affirmative. First, we
will apply the central limit theorem for stationary, strongly mixing
sequences, and this will allow us to evaluate the variance parameter.
Then, using a simple coupling argument, we will treat the general case
in which the initial state is fixed but arbitrary.

It will be convenient to index time by ${\bf Z}$, the set of integers.
So we let $\{R_n\}_{n\in\mathbf{Z}}$ be independent, uniformly bounded,
nonnegative random variables, with $\{R_n\dvtx  n-1\equiv i$ (mod $I)\}$
identically distributed as the payout distribution in cam position
$i\in\{0,1,\ldots,I-1\}$. We interpret $\{R_n\}_{n\in\mathbf{Z}}$
as the
sequence of payouts of the slot machine \textit{excluding} the
Futurity awards. Thus,
\[
\mathrm{P}(R_n>0)=p_{n-1},\qquad  \mathrm{E}[R_n]=\mu_{n-1},\qquad
\operatorname{Var}(R_n)=\sigma_{n-1}^2
\]
for all $n\in\mathbf{Z}$, provided we extend these parameters periodically;
for example, $p_{i+mI}:=p_i$ for all $i\in\{0,1,\ldots,I-1\}$ and
$m\in\mathbf{Z}$.

Notice that we can define the Markov chain $\{(X_n,Y_n)\}_{n\in\mathbf{Z}}$
as a nonrandom function of $\{R_n\}_{n\in\mathbf{Z}}$. Indeed,
$X_n=i\in\{
0,1,\ldots,I-1\}$ if $n\equiv i$ (mod $I$), so $\{X_n\}_{n\in\mathbf
{Z}}$ is
deterministic, and $Y_n=j\in\{0,1,\ldots,J-1\}$ if
\[
R_{n-kJ-j}>0,\quad   R_{n-kJ-j+1}=\cdots=R_n=0\qquad  \mbox{for some }k\ge0.
\]
To take the Futurity awards into account, we define $\{R_n^*\}_{n\in
\mathbf{Z}
}$ by
%
\begin{eqnarray}\label{R_n^*}
R_n^*&:=&R_n+J\cdot1_{\{Y_{n-1}=J-1,  R_n=0\}}\nonumber\\
& =&R_n+J\sum_{k=1}^\infty1_{\{R_{n-kJ}>0,R_{n-kJ+1}=\cdots
=R_{n-1}=R_n=0\}}\\
& =&u(\ldots,R_{n-2},R_{n-1},R_n), \qquad  n\in\mathbf{Z}, \nonumber
\end{eqnarray}
for some nonrandom function $u$.

The sequence $\{R_n\}_{n\in\mathbf{Z}}$ is independent but not identically
distributed, so we consider the sequence of random vectors
\[
\mathbf{R}_k:=\bigl(R_{kI+1},\ldots,R_{(k+1)I}\bigr),\qquad   k\in\mathbf{Z},
\]
which \textit{is} i.i.d., hence by (\ref{R_n^*}),
\[
\mathbf{R}_k^*:=\bigl(R_{kI+1}^*,\ldots,R_{(k+1)I}^*\bigr), \qquad  k\in\mathbf{Z},
\]
is a stationary sequence. In particular, the sequence
\[
S_k:=R_{kI+1}+\cdots+R_{(k+1)I}, \qquad  k\in\mathbf{Z},
\]
is also i.i.d., and the sequence
\[
S_k^*:=R_{kI+1}^*+\cdots+R_{(k+1)I}^*,\qquad   k\in\mathbf{Z},
\]
is also stationary, despite the fact that the Markov chain $\{
(X_n,Y_n)\}_{n\in\mathbf{Z}}$ is not stationary in this construction. $S_k$
and $S_k^*$ represent the total payout, excluding and including the
Futurity awards, respectively, over the segment of $I$ consecutive
coups numbered $kI+1,\ldots,(k+1)I$.

We claim that the stationary sequence $\{S_k^*\}_{k\in\mathbf{Z}}$ is
strongly mixing, that is, the quantities
%
\begin{equation}\label{alpha}
\alpha(m):=\sup_{A\in\sigma(S_k^*:k\le-m), B\in\sigma
(S_k^*:k\ge0)}|\mathrm{P}(A\cap B)-\mathrm{P}(A)\mathrm{P}(B)|
\end{equation}
satisfy $\alpha(m)\to0$ as $m\to\infty$. For $m\ge2$, let $C_m:=\{
R_k>0$ for some $k\in\{-(m-1)I+1,-(m-1)I+2,\ldots,0\}$. Then, with
$A$ and $B$ as in (\ref{alpha}), $A$ is independent of $B\cap C_m$, so
\begin{eqnarray*}
&&|\mathrm{P}(A\cap B)-\mathrm{P}(A)\mathrm{P}(B)|\nonumber\\
&&\qquad \le|\mathrm{P}(A\cap B\cap C_m)-\mathrm{P}(A)\mathrm{P}(B\cap
C_m)|\nonumber\\
&&\qquad \quad   {}+|\mathrm{P}(A\cap B\cap C_m^c)-\mathrm
{P}(A)\mathrm{P}(B\cap
C_m^c)|\nonumber\\
&&\qquad  =|\mathrm{P}(A\cap B\cap C_m^c)-\mathrm{P}(A)\mathrm{P}(B\cap
C_m^c)|\nonumber\\
&& \qquad\le\mathrm{P}(C_m^c)\\
&&\qquad =\mathrm
{P}\bigl(R_{-(m-1)I+1}=R_{-(m-1)I+2}=\cdots=R_{0}=0\bigr)\\
&&\qquad =Q^{m-1},
\end{eqnarray*}
and this shows that $\alpha(m)$ converges to 0 geometrically fast.

Letting $\bar\mu:=\mathrm{E}[S_0^*]$ and noting that the random
variables of
interest are uniformly bounded, the central limit theorem for
stationary, strongly mixing sequences [e.g., Bradley (\citeyear{B}), Theorem
10.3] tells us that
\[
{S_0^*+\cdots+S_{m-1}^*-m\bar\mu\over\sqrt{m\bar\sigma
^2}}\stackrel{d}{\to} N(0,1),
\]
provided
\[
\bar\sigma^2:=\operatorname{Var}(S_0^*)+2\sum_{m=1}^\infty
\operatorname{Cov}(S_0^*, S_m^*)>0.
\]
We now evaluate $\bar\sigma^2$.

First, we will frequently encounter
\begin{eqnarray*}
&&\mathrm{P}(Y_{i-1}=J-1, R_i=0)\\
&&\qquad  {}=\sum_{k=1}^\infty\mathrm{P}(R_{i-kJ}>0,
R_{i-kJ+1}=\cdots
=R_i=0)\nonumber\\
&& \qquad {}={1\over1-Q}\sum_{k=1}^d q_{i-1}\cdots q_{i-kJ}p_{i-kJ-1}\\
&& \qquad {}=:P_{i-1}
\end{eqnarray*}
for $i=1,2,\ldots,I$. For example,
\begin{eqnarray*}
\bar\mu&:=&\mathrm{E}[S_0^*]=\mathrm{E}[R_1^*+\cdots
+R_I^*]\nonumber\\
& =&\sum_{i=1}^{I}\mathrm{E}[R_i^*]=\sum_{i=1}^{I}\mathrm
{E}\bigl[R_i+J\cdot1_{\{
Y_{i-1}=J-1, R_i=0\}}\bigr]\\
& =&\sum_{i=1}^{I}(\mu_{i-1}+JP_{i-1})=\sum_{i=0}^{I-1}\mu_i+J\sum
_{i=0}^{I-1}P_i.
\end{eqnarray*}
Next, for $i=1,2,\ldots,I$,
\begin{eqnarray*}
\operatorname{Var}(R_i^*)&=&\operatorname{Var}\bigl(R_i+J\cdot1_{\{
Y_{i-1}=J-1, R_i=0\}}\bigr)\nonumber\\
&=&\operatorname{Var}(R_i)+2J\operatorname{Cov}\bigl(R_i,1_{\{
Y_{i-1}=J-1, R_i=0\}}\bigr)\nonumber\\
&&  {}+J^2\operatorname{Var}\bigl(1_{\{Y_{i-1}=J-1, R_i=0\}
}\bigr)\nonumber\\
&=&\operatorname{Var}(R_i)-2J\mathrm{E}[R_i]\mathrm{P}(Y_{i-1}=J-1,
R_i=0)\nonumber\\
&& {}+J^2\mathrm{P}(Y_{i-1}=J-1, R_i=0)\bigl(1-\mathrm
{P}(Y_{i-1}=J-1,
R_i=0)\bigr)\nonumber\\
&=&\sigma_{i-1}^2-2J\mu_{i-1}P_{i-1}+J^2P_{i-1}(1-P_{i-1}),
\end{eqnarray*}
and, for $1\le i<j\le I$,
\begin{eqnarray*}
\operatorname{Cov}(R_i^*,R_j^*)&=&\operatorname{Cov}\bigl(R_i+J\cdot1_{\{
Y_{i-1}=J-1, R_i=0\}
},R_j+J\cdot1_{\{Y_{j-1}=J-1, R_j=0\}}\bigr)\nonumber\\
&=&J\operatorname{Cov}\bigl(R_i,1_{\{Y_{j-1}=J-1, R_j=0\}}\bigr)\nonumber\\
&& {}+J^2\operatorname{Cov}\bigl(1_{\{Y_{i-1}=J-1, R_i=0\}},1_{\{
Y_{j-1}=J-1,
R_j=0\}}\bigr)\nonumber\\
&=&J\bigl\{\mathrm{E}\bigl[R_i 1_{\{Y_{j-1}=J-1, R_j=0\}}\bigr]-\mathrm
{E}[R_i]\mathrm{P}(Y_{j-1}=J-1,
R_j=0)\bigr\}\nonumber\\
&& {}+J^2[\mathrm{P}(Y_{i-1}=J-1, R_i=0, Y_{j-1}=J-1,
R_j=0)\nonumber
\\
&& \hspace*{25pt} {}-\mathrm{P}(Y_{i-1}=J-1, R_i=0)\mathrm
{P}(Y_{j-1}=J-1,
R_j=0)]\nonumber\\
&=&J\mu_{i-1}\biggl(\sum_{1\le k<(j-i)/J}q_{j-1}\cdots
q_{j-kJ}p_{j-kJ-1}\nonumber\\
&& \hspace*{82pt}   {}+\Delta_{ij}q_{j-1}\cdots
q_i-P_{j-1}\biggr)\nonumber\\
&& {}+J^2\biggl(\sum_{1\le k<(j-i)/J}q_{j-1}\cdots
q_{j-kJ}p_{j-kJ-1}P_{i-1}\nonumber\\
&&  \hspace*{31pt}{}+\Delta_{ij}\sum_{k>(j-i)/J}q_{j-1}\cdots
q_{j-kJ}p_{j-kJ-1}-P_{i-1}P_{j-1}\biggr)\nonumber\\
&=:&A_{ij},
\end{eqnarray*}
where $\Delta_{ij}:=1$ if $j-i\equiv0$ (mod $J$) and $:=0$ otherwise,
and the infinite series in the definition of $A_{ij}$ can be expressed
as the finite sum
\[
\sum_{(j-i)/J<k\le d}q_{j-1}\cdots q_{j-kJ}p_{j-kJ-1}+Q P_{j-1}
\]
when $j-i\equiv0$ (mod $J$). We conclude that
\begin{eqnarray}\label{var}
   \operatorname{Var}(S_0^*)&=&\operatorname
{Var}(R_1^*+\cdots+R_I^*)=\sum_{i=1}^I\operatorname{Var}
(R_i^*)+2\mathop{\sum\sum}_{1\le i<j\le I}\operatorname
{Cov}(R_i^*,R_j^*)\nonumber\\[-8pt]\\[-8pt]
&=&\sum_{i=0}^{I-1}[\sigma_i^2-2J\mu_iP_i+J^2P_i(1-P_i)]+2\mathop
{\sum\sum}_{1\le i<j\le I}A_{ij}.\nonumber
\end{eqnarray}
Notice that this formula depends solely on the basic parameters ($I$,
$J$, $p_i$, $\mu_i$ and~$\sigma_i^2$).

Next, for $i,j=1,2,\ldots,I$ and $m\ge1$,
\begin{eqnarray*}
&&\operatorname{Cov}(R_i^*,R_{mI+j}^*)\nonumber\\
&& \qquad =\operatorname{Cov}\bigl(R_i+J\cdot1_{\{Y_{i-1}=J-1, R_i=0\}},
R_{mI+j}+J\cdot1_{\{Y_{mI+j-1}=J-1, R_{mI+j}=0\}}\bigr)\nonumber\\
&& \qquad =J \operatorname{Cov}\bigl(R_i, 1_{\{Y_{mI+j-1}=J-1,
R_{mI+j}=0\}}\bigr)\nonumber
\\
&&\quad \qquad   {}+J^2 \operatorname{Cov}\bigl(1_{\{Y_{i-1}=J-1, R_i=0\}
}, 1_{\{
Y_{mI+j-1}=J-1, R_{mI+j}=0\}}\bigr).
\end{eqnarray*}
Now
\begin{eqnarray*}
&&\operatorname{Cov}\bigl(R_i, 1_{\{Y_{mI+j-1}=J-1, R_{mI+j}=0\}
}\bigr)\nonumber\\
&&\qquad =\sum_{1\le k\le md+(j-i)/J}\mathrm{E}\bigl[R_i 1_{\{
R_{mI-kJ+j}>0,
R_{mI-kJ+j+1}=\cdots=R_{mI+j}=0\}}\bigr]\nonumber\\
&& \quad \qquad     {}-\mu_{i-1}P_{j-1}\nonumber\\
&& \qquad =\sum_{1\le k<md+(j-i)/J}q_{j-1}\cdots q_{j-kJ}p_{j-kJ-1}\mu
_{i-1}\nonumber\\
&& \quad \qquad   {}+\Delta_{ij}q_{mI+j-1}\cdots q_i\mu
_{i-1}-\mu_{i-1}P_{j-1}\nonumber\\
&&\qquad =\mu_{i-1}\Biggl(\sum_{k=1}^d q_{j-1}\cdots
q_{j-kJ}p_{j-kJ-1}(1+Q+\cdots+Q^{m-2})\nonumber\\
&&\hspace*{60pt}   {}+Q^{m-1}\sum_{1\le k<d+(j-i)/J}q_{j-1}\cdots
q_{j-kJ}p_{j-kJ-1}\nonumber\\
&& \hspace*{97pt}  {}+\Delta_{ij}q_{j-1}\cdots q_0
Q^{m-1}q_{I-1}\cdots q_i-P_{j-1}\Biggr)\nonumber\\
&& \qquad =\mu_{i-1}\biggl( -P_{j-1}+\sum_{1\le
k<d+(j-i)/J}q_{j-1}\cdots q_{j-kJ}p_{j-kJ-1}\nonumber\\
&&   \hspace*{152pt}  {}+\Delta_{ij}q_{I-1}\cdots q_i
q_{j-1}\cdots q_0\biggr)Q^{m-1}\nonumber\\
&&\qquad =:B_{ij}Q^{m-1},
\end{eqnarray*}
where $q_{I-1}\cdots q_i:=1$ if $i=I$, and
\begin{eqnarray*}
&&\operatorname{Cov}\bigl(1_{\{Y_{i-1}=J-1, R_i=0\}}, 1_{\{
Y_{mI+j-1}=J-1, R_{mI+j}=0\}
}\bigr)\nonumber\\
&& \qquad =\mathrm{P}(Y_{i-1}=J-1, R_i=0, Y_{mI+j-1}=J-1,
R_{mI+j}=0)-P_{i-1}P_{j-1}\nonumber\\
&&  \qquad =\sum_{1\le k<md+(j-i)/J}  P_{i-1}\mathrm
{P}(R_{mI-kJ+j}>0,
R_{mI-kJ+j+1}=\cdots=R_{mI+j}=0)\nonumber\\
&& \quad \qquad    {}+\Delta_{ij}\sum_{k>md+(j-i)/J} \mathrm
{P}(R_{mI-kJ+j}>0,
R_{mI-kJ+j+1}=\cdots=R_{mI+j}=0)\nonumber\\
&& \quad \qquad       {}-P_{i-1}P_{j-1}\nonumber\\
&&  \qquad =\sum_{1\le k<md+(j-i)/J}q_{j-1}\cdots
q_{j-kJ}p_{j-kJ-1}P_{i-1}\nonumber\\
&& \quad \qquad     {}+\Delta_{ij}\sum_{k>md+(j-i)/J}q_{j-1}\cdots
q_{j-kJ}p_{j-kJ-1}-P_{i-1}P_{j-1}\\
&& \qquad =\sum_{k=1}^d q_{j-1}\cdots q_{j-kJ}p_{j-kJ-1}(1+Q+\cdots
+Q^{m-2})P_{i-1}\nonumber\\
&&  \quad \qquad    {}+Q^{m-1}\sum_{1\le
k<d+(j-i)/J}q_{j-1}\cdots q_{j-kJ}p_{j-kJ-1}P_{i-1}\nonumber\\
&& \quad \qquad   {}+\Delta_{ij}Q^{m-1}\sum_{k>d+(j-i)/J}q_{j-1}\cdots
q_{j-kJ}p_{j-kJ-1}-P_{i-1}P_{j-1}\\
&&  \qquad =\biggl( -P_{i-1}P_{j-1}+\sum_{1\le k<d+(j-i)/J}q_{j-1}\cdots
q_{j-kJ}p_{j-kJ-1}P_{i-1}\nonumber\\
&&  \hspace*{106pt}     {}+\Delta_{ij}\sum
_{k>d+(j-i)/J}q_{j-1}\cdots q_{j-kJ}p_{j-kJ-1}\biggr)Q^{m-1}\nonumber\\
&&  \qquad =:C_{ij}Q^{m-1},
\end{eqnarray*}
and the infinite series in the definition of $C_{ij}$ can be expressed as
\[
Q\sum_{(j-i)/J<k\le d}q_{j-1}\cdots q_{j-kJ}p_{j-kJ-1}+Q^2 P_{j-1}
\]
when $j-i\equiv0$ (mod $J$) and $j\ge i$, and as
\[
\sum_{d+(j-i)/J<k\le d}q_{j-1}\cdots q_{j-kJ}p_{j-kJ-1}+Q P_{j-1}
\]
when $j-i\equiv0$ (mod $J$) and $j<i$. We conclude that
\begin{eqnarray*}
\operatorname{Cov}(S_0^*,S_m^*)&=&\operatorname{Cov}\bigl(R_1^*+\cdots
+R_I^*,R_{mI+1}^*+\cdots
+R_{(m+1)I}^*\bigr)\nonumber\\
&=&\sum_{i=1}^I\sum_{j=1}^I\operatorname{Cov}(R_i^*,R_{mI+j}^*)\\
&=&\sum
_{i=1}^I\sum
_{j=1}^I (JB_{ij}Q^{m-1}+J^2C_{ij}Q^{m-1}),
\end{eqnarray*}
and hence that
%
\begin{eqnarray}\label{covsum}
\sum_{m=1}^\infty\operatorname{Cov}(S_0^*,S_m^*)&=&{J\over1-Q}\sum
_{i=1}^I\sum
_{j=1}^I(B_{ij}+JC_{ij}).
\end{eqnarray}
Again, this formula depends solely on the basic parameters. Summing
(\ref{var}) and twice (\ref{covsum}), we obtain $\bar\sigma^2$.

Finally, we observe that the central limit theorem
for the stationary sequence $\{S_k^*\}_{k\in\mathbf{Z}}$ yields a central
limit theorem for $\{R_n^*\}_{n\in\mathbf{Z}}$ as well. Indeed, with
\[
\mu^*=\bar\mu/I \quad \mbox{and}\quad  (\sigma^*)^2:=\bar\sigma^2/I,
\]
we find that
%
\begin{equation}\label{couple}
{R_1^*+\cdots+R_n^*-n\mu^*\over\sqrt{n(\sigma^*)^2}}-{S_0^*+\cdots
+S_{\lfloor n/I\rfloor-1}^*-\lfloor n/I\rfloor\bar\mu\over\sqrt
{(n/I)\bar\sigma^2}}
\end{equation}
tends to 0 a.s.\ as $n\to\infty$ because the difference between the
numerators, namely
\[
R_{I\lfloor n/I\rfloor+1}^*+\cdots+R_n^*-(n-I\lfloor n/I\rfloor)\mu^*,
\]
is uniformly bounded in $n$ and the denominators are equal. Thus,
\[
{R_1^*+\cdots+R_n^*-n\mu^*\over\sqrt{n(\sigma^*)^2}}\stackrel
{d}{\to} N(0,1).
\]

We can go one step further and derive a central limit theorem for $\{
\hat R_n^*\}_{n\ge0}$ with $(\hat X_0,\hat Y_0)=(i_0,j_0)$ specified,
where the hats on $\hat R_n^*$, $\hat X_0$ and $\hat Y_0$ distinguish
them from the $R_n^*$, $X_0$ and $Y_0$ already defined. The idea of the
proof is the same as in~(\ref{couple}). We define $\hat
R_n:=R_{n+i_0}$ for $n\ge1$, and we define $(\hat X_n,\hat Y_n)$ for
$n\ge1$ in terms of $(\hat X_0,\hat Y_0)=(i_0,j_0)$ and $1_{\{\hat
R_1=0\}},\ldots,1_{\{\hat R_n=0\}}$ in the usual way. Then
\[
{\hat R_1^*+\cdots+\hat R_n^*-n\mu^*\over\sqrt{n(\sigma
^*)^2}}-{R_1^*+\cdots+R_n^*-n\mu^*\over\sqrt{n(\sigma^*)^2}}
\]
tends to 0 a.s.\ as $n\to\infty$ because $\hat R_n^*=R_{n+i_0}^*$
unless $\hat R_1^*=\cdots=\hat R_{n-1}^*=0$. In\vspace*{-1pt} words, the sequences
$\hat R_1^*,\ldots,\hat R_n^*$ and $R_1^*,\ldots,R_n^*$ differ only
by a shift (of $i_0$ terms), once the $\hat Y$ process and the shifted
$Y$ process couple, which occurs after the first win. We have therefore
established the following central limit theorem.

\begin{theorem}\label{teo3}
Let $\hat R_1^*,\hat R_2^*,\ldots$ be the sequence of payouts of the
generalized Futurity slot machine starting in an arbitrary initial
state. Then
\[
{\hat R_1^*+\cdots+\hat R_n^*-n\mu^*\over\sqrt{n(\sigma
^*)^2}}\stackrel{d}{\to} N(0,1).
\]
\end{theorem}

\section{Numerical results for the Futurity}\label{sec5}

The Futurity was in production from 1936 to 1941. (After December 7, 1941,
Mills Novelty stopped producing slot machines and became a defense
contractor for the duration of the war. When it resumed slot production
in 1945, it did so with new designs.) In particular, there were minor
variations in the payouts and reel strip labels used with the machine,
but the fundamental properties, the Futurity award and the periodic
mode changes, are common to every Mills Futurity. The precise version
we consider here is the one described by Geddes (\citeyear{G}).

To simplify matters, we code the six symbols as lemon${}=0$,
cherry${}=1$, orange${}=2$, plum${}=3$, bell${}=4$ and bar${}=5$. The
pay table can then be described by the function $p\dvtx \{0,1,2,3,4,5\}
^3\mapsto{\bf Z}_+$ given by $p(5,5,5):=150$, $p(4,4,4)=p(4,4,5):=18$,
$p(3,3,3)=p(3,3,5):=14$, $p(2,2,2)=p(2,2,5):=10$,
$p(1,1,0)=p(1,1,4):=5$ and $p(1,1,2)=p(1,1,3)=p(1,1,5):=3$; otherwise
$p:=0$. The three reel strips can be described as follows, in which the
symbols in odd-numbered positions are italicized for convenience:
\begin{eqnarray*}
&&\mbox{reel 1:}\quad  \mathit{1},5,\mathit{1},2,\mathit{1},5,\mathit
{1},5,\mathit{1},3,\mathit{1},2,\mathit{5},1,\mathit{4},3,\mathit
{1},5,\mathit{1},2,\nonumber\\
&&\mbox{reel 2:}\quad  \mathit{1},4,\mathit{1},3,\mathit{1},4,\mathit
{1},2,\mathit{1},4,\mathit{1},4,\mathit{1},2,\mathit{1},2,\mathit
{4},1,\mathit{5},4,\nonumber\\
&&\mbox{reel 3:}\quad  \mathit{3},4,\mathit{2},0,\mathit{3},4,\mathit
{2},0,\mathit{4},0,\mathit{2},3,\mathit{2},4,\mathit{2},4,\mathit
{5},2,\mathit{3},5.\nonumber
\end{eqnarray*}
Table \ref{Table1} summarizes the relevant information from these reel strips.
Of course, the reels operate independently, and the 10 possible
positions at which each reel can stop (given the mode) are assumed
equally likely.

\begin{table}
\caption{Reel strip inventories for the Futurity in
both modes}\label{Table1}
\begin{tabular*}{\textwidth}{@{\extracolsep{\fill}}lcccccc@{}}
\hline
& \multicolumn{3}{c}{\textbf{Mode E}} & \multicolumn{3}{c@{}}{\textbf{Mode O}} \\[-6pt]
& \multicolumn{3}{@{\hspace*{-57pt}}c@{\hspace*{-57pt}}}{\hrulefill} & \multicolumn{3}{@{\hspace*{-56pt}}c@{\hspace*{-63pt}}}{\hrulefill} \\
\textbf{Symbol} & \textbf{Reel 1} & \textbf{Reel 2} & \textbf{Reel 3} & \textbf{Reel 1} & \textbf{Reel 2} & \textbf{Reel 3} \\
\hline
Lemon $(=0)$ & \phantom{0}0 & \phantom{0}0 & \phantom{0}3 & \phantom{0}0 & \phantom{0}0 & \phantom{0}0 \\
Cherry $(=1)$ & \phantom{0}1 & \phantom{0}1 & \phantom{0}0 & \phantom{0}8 & \phantom{0}8 & \phantom{0}0 \\
Orange $(=2)$ & \phantom{0}3 & \phantom{0}3 & \phantom{0}1 & \phantom{0}0 & \phantom{0}0 & \phantom{0}5 \\
Plum $(=3)$ & \phantom{0}2 & \phantom{0}1 & \phantom{0}1 & \phantom{0}0 & \phantom{0}0 & \phantom{0}3 \\
Bell $(=4)$ & \phantom{0}0 & \phantom{0}5 & \phantom{0}4 & \phantom{0}1 & \phantom{0}1 & \phantom{0}1 \\
Bar $(=5)$ & \phantom{0}4 & \phantom{0}0 & \phantom{0}1 & \phantom{0}1 & \phantom{0}1 & \phantom{0}1 \\
[3pt]
Total & 10 & 10 & 10 & 10 & 10 & 10 \\
\hline
\end{tabular*}
\end{table}

With $f_\mathrm{E}(i,j)$ denoting the frequency of symbol $i$ on reel
$j$ in
mode E (see Table~\ref{Table1}), we find that the mean payout in mode
E is
\[
\mu_\mathrm{E}={1\over(10)^3}\sum_{i_1=0}^5\sum_{i_2=0}^5\sum
_{i_3=0}^5f_\mathrm{E}(i_1,1)f_\mathrm{E}(i_2,2)f_\mathrm
{E}(i_3,3)p(i_1,i_2,i_3)=0.28.
\]
Similarly, the mean payout in mode O is $\mu_\mathrm{O}=2.234$. Certainly,
these numbers justify our descriptions of mode E as ``tight'' and mode
O as ``loose,''
as do the facts that the probability of a nonzero payout in mode E,
other than a Futurity award, is $p_\mathrm{E}=0.032$, and the corresponding
probability in mode O is $p_\mathrm{O}=0.643$. See Table \ref{Table2}.

\begin{table}[b]
\caption{Payout frequencies and statistics for the
Futurity in both modes, excluding Futurity awards. Results~are exact
(no rounding)}\label{Table2}
\begin{tabular*}{\textwidth}{@{\extracolsep{\fill}}lcc@{}}
\hline
\textbf{Payout} & \textbf{Mode E} & \textbf{Mode O} \\
\hline
\phantom{00}0 & \phantom{0}968 & \phantom{0}357 \\
\phantom{00}3 & \phantom{000}3 & \phantom{0}576 \\
\phantom{00}5 & \phantom{000}7 & \phantom{00}64 \\
\phantom{0}10 & \phantom{00}18 & \phantom{000}0 \\
\phantom{0}14 & \phantom{000}4 & \phantom{000}0 \\
\phantom{0}18 & \phantom{000}0 & \phantom{000}2 \\
150 & \phantom{000}0 & \phantom{000}1 \\
[3pt]
Total & 1000 & 1000 \\
[6pt]
Mean payout & $\mu_\mathrm{E}=0.28$ & $\mu_\mathrm{O}=2.234$ \\
Variance of payout & $\sigma_\mathrm{E}^2=2.7076$ & $\sigma_\mathrm{O}^2=24.941244$\\
Probability of nonzero payout & $p_\mathrm{E}=0.032$ & $p_\mathrm{O}=0.643$ \\
\hline
\end{tabular*}
\end{table}

\begin{sidewaystable}
\tablewidth=\textheight
\tablewidth=\textwidth
\caption{Stationary distribution of the Markov chain,
rounded to six decimal places. Rows indicate cam position, and columns
indicate pointer position. Entries~greater than 1$/$100 are shaded}\label{Table3}
\begin{tabular*}{\textwidth}{@{\extracolsep{\fill}}l
>{\columncolor[gray]{0.85}}c
>{\columncolor[gray]{0.85}}c
>{\columncolor[gray]{0.85}}c
>{\columncolor[gray]{0.85}}c
>{\columncolor[gray]{0.85}}c
>{\columncolor[gray]{0.85}}c
>{\columncolor[gray]{0.85}}c
>{\columncolor[gray]{0.85}}c
>{\columncolor[gray]{0.85}}c
>{\columncolor[gray]{0.85}}c
c@{}}
\hline
& \multicolumn{1}{>{\columncolor{white}}c}{\textbf{0}} & \multicolumn{1}{>{\columncolor{white}}c}{\textbf{1}}
 & \multicolumn{1}{>{\columncolor{white}}c}{\textbf{2}} & \multicolumn{1}{>{\columncolor{white}}c}{\textbf{3}}
  & \multicolumn{1}{>{\columncolor{white}}c}{\textbf{4}} & \multicolumn{1}{>{\columncolor{white}}c}{\textbf{5}}
   & \multicolumn{1}{>{\columncolor{white}}c}{\textbf{6}} & \multicolumn{1}{>{\columncolor{white}}c}{\textbf{7}}
    & \multicolumn{1}{>{\columncolor{white}}c}{\textbf{8}} & \multicolumn{1}{>{\columncolor{white}}c}{\textbf{9}} & \textbf{Sum} \\
\hline
0 &0.071306 &\multicolumn{1}{>{\columncolor{white}}c}{0.001267} &\multicolumn{1}{>{\columncolor{white}}c}{0.001226}
 &\multicolumn{1}{>{\columncolor{white}}c}{0.001187}
&0.023090 &\multicolumn{1}{>{\columncolor{white}}c}{0.000410} &\multicolumn{1}{>{\columncolor{white}}c}{0.000397}
 &\multicolumn{1}{>{\columncolor{white}}c}{0.000384}
 &\multicolumn{1}{>{\columncolor{white}}c}{0.000372}
&\multicolumn{1}{>{\columncolor{white}}c}{0.000360} & $1/10$ \\
1 &\multicolumn{1}{>{\columncolor{white}}c}{0.003549} &0.069024 &\multicolumn{1}{>{\columncolor{white}}c}{0.001226}
&\multicolumn{1}{>{\columncolor{white}}c}{0.001187}
&\multicolumn{1}{>{\columncolor{white}}c}{0.001149}
&0.022351 &\multicolumn{1}{>{\columncolor{white}}c}{0.000397} &\multicolumn{1}{>{\columncolor{white}}c}{0.000384}
 &\multicolumn{1}{>{\columncolor{white}}c}{0.000372} &\multicolumn{1}{>{\columncolor{white}}c}{0.000360} &
$1/10$ \\
2 &\multicolumn{1}{>{\columncolor{white}}c}{0.003549}
 &\multicolumn{1}{>{\columncolor{white}}c}{0.003435}
  &0.066815 &\multicolumn{1}{>{\columncolor{white}}c}{0.001187}
   &\multicolumn{1}{>{\columncolor{white}}c}{0.001149}
&\multicolumn{1}{>{\columncolor{white}}c}{0.001112} &0.021636 &\multicolumn{1}{>{\columncolor{white}}c}{0.000384}
 &\multicolumn{1}{>{\columncolor{white}}c}{0.000372} &\multicolumn{1}{>{\columncolor{white}}c}{0.000360} &
$1/10$ \\
3 &\multicolumn{1}{>{\columncolor{white}}c}{0.003549} &\multicolumn{1}{>{\columncolor{white}}c}{0.003435}
 &\multicolumn{1}{>{\columncolor{white}}c}{0.003325} &0.064677
 &\multicolumn{1}{>{\columncolor{white}}c}{0.001149}
&\multicolumn{1}{>{\columncolor{white}}c}{0.001112} &\multicolumn{1}{>{\columncolor{white}}c}{0.001077} &0.020943
 &\multicolumn{1}{>{\columncolor{white}}c}{0.000372} &\multicolumn{1}{>{\columncolor{white}}c}{0.000360} &
$1/10$ \\
4 &\multicolumn{1}{>{\columncolor{white}}c}{0.003549} &\multicolumn{1}{>{\columncolor{white}}c}{0.003435}
 &\multicolumn{1}{>{\columncolor{white}}c}{0.003325} &\multicolumn{1}{>{\columncolor{white}}c}{0.003219} &0.062608
&\multicolumn{1}{>{\columncolor{white}}c}{0.001112} &\multicolumn{1}{>{\columncolor{white}}c}{0.001077}
 &\multicolumn{1}{>{\columncolor{white}}c}{0.001042} &0.020273 &\multicolumn{1}{>{\columncolor{white}}c}{0.000360} &
$1/10$ \\
5 &\multicolumn{1}{>{\columncolor{white}}c}{0.003549} &\multicolumn{1}{>{\columncolor{white}}c}{0.003435}
 &\multicolumn{1}{>{\columncolor{white}}c}{0.003325} &\multicolumn{1}{>{\columncolor{white}}c}{0.003219}
  &\multicolumn{1}{>{\columncolor{white}}c}{0.003116} &0.060604 &\multicolumn{1}{>{\columncolor{white}}c}{0.001077}
   &\multicolumn{1}{>{\columncolor{white}}c}{0.001042} &\multicolumn{1}{>{\columncolor{white}}c}{0.001009}
    &0.019624 & $1/10$ \\
6 &0.071306 &\multicolumn{1}{>{\columncolor{white}}c}{0.001267} &\multicolumn{1}{>{\columncolor{white}}c}{0.001226}
 &\multicolumn{1}{>{\columncolor{white}}c}{0.001187}
 &\multicolumn{1}{>{\columncolor{white}}c}{0.001149}
&\multicolumn{1}{>{\columncolor{white}}c}{0.001112} &0.021636 &\multicolumn{1}{>{\columncolor{white}}c}{0.000384}
 &\multicolumn{1}{>{\columncolor{white}}c}{0.000372} &\multicolumn{1}{>{\columncolor{white}}c}{0.000360} &
$1/10$ \\
7 &\multicolumn{1}{>{\columncolor{white}}c}{0.003549} &0.069024 &\multicolumn{1}{>{\columncolor{white}}c}{0.001226}
 &\multicolumn{1}{>{\columncolor{white}}c}{0.001187}
 &\multicolumn{1}{>{\columncolor{white}}c}{0.001149}
&\multicolumn{1}{>{\columncolor{white}}c}{0.001112} &\multicolumn{1}{>{\columncolor{white}}c}{0.001077} &0.020943
 &\multicolumn{1}{>{\columncolor{white}}c}{0.000372} &\multicolumn{1}{>{\columncolor{white}}c}{0.000360} &
$1/10$ \\
8 &\multicolumn{1}{>{\columncolor{white}}c}{0.003549} &\multicolumn{1}{>{\columncolor{white}}c}{0.003435} &0.066815
 &\multicolumn{1}{>{\columncolor{white}}c}{0.001187}
 &\multicolumn{1}{>{\columncolor{white}}c}{0.001149}
&\multicolumn{1}{>{\columncolor{white}}c}{0.001112} &\multicolumn{1}{>{\columncolor{white}}c}{0.001077}
 &\multicolumn{1}{>{\columncolor{white}}c}{0.001042} &0.020273
  &\multicolumn{1}{>{\columncolor{white}}c}{0.000360} &
$1/10$ \\
9 &\multicolumn{1}{>{\columncolor{white}}c}{0.003549} &\multicolumn{1}{>{\columncolor{white}}c}{0.003435}
 &\multicolumn{1}{>{\columncolor{white}}c}{0.003325} &0.064677
 &\multicolumn{1}{>{\columncolor{white}}c}{0.001149}
&\multicolumn{1}{>{\columncolor{white}}c}{0.001112} &\multicolumn{1}{>{\columncolor{white}}c}{0.001077}
 &\multicolumn{1}{>{\columncolor{white}}c}{0.001042} &\multicolumn{1}{>{\columncolor{white}}c}{0.001009} &0.019624 &
$1/10$ \\
[3pt]
Sum &\multicolumn{1}{>{\columncolor{white}}c}{0.171001} &\multicolumn{1}{>{\columncolor{white}}c}{0.161193}
 &\multicolumn{1}{>{\columncolor{white}}c}{0.151837} &\multicolumn{1}{>{\columncolor{white}}c}{0.142915}
  &\multicolumn{1}{>{\columncolor{white}}c}{0.096857} &\multicolumn{1}{>{\columncolor{white}}c}{0.091152}
   &\multicolumn{1}{>{\columncolor{white}}c}{0.050526}
&\multicolumn{1}{>{\columncolor{white}}c}{0.047593} &\multicolumn{1}{>{\columncolor{white}}c}{0.044797}
 &\multicolumn{1}{>{\columncolor{white}}c}{0.042130} & \\
\hline
\end{tabular*}
\end{sidewaystable}

With the statistics of Table \ref{Table2}, we can define
\begin{eqnarray*}
(p_0,p_1,\ldots,p_9)&:=&(p_\mathrm{E},p_\mathrm{E},p_\mathrm
{E},p_\mathrm{E},p_\mathrm{E},p_\mathrm{O},p_\mathrm{E}
,p_\mathrm{E},p_\mathrm{E},p_\mathrm{O}),\\
(\mu_0,\mu_1,\ldots,\mu_{9})&:=&(\mu_\mathrm{E},\mu_\mathrm
{E},\mu_\mathrm{E},\mu_\mathrm{E}
,\mu_\mathrm{E},\mu_\mathrm{O},\mu_\mathrm{E},\mu_\mathrm
{E},\mu_\mathrm{E},\mu_\mathrm{O}),\\
(\sigma_0^2,\sigma_1^2,\ldots,\sigma_{9}^2)&:=&(\sigma_\mathrm
{E}^2,\sigma
_\mathrm{E}^2,\sigma_\mathrm{E}^2,\sigma_\mathrm{E}^2,\sigma
_\mathrm{E}^2,\sigma_\mathrm{O}^2,\sigma
_\mathrm{E}^2,\sigma_\mathrm{E}^2,\sigma_\mathrm{E}^2,\sigma
_\mathrm{O}^2),
\end{eqnarray*}
and $q_i:=1-p_i$ for $i=0,1,\ldots,9$. With $I=J=10$ (in particular,
$d=1$), we can apply Theorem \ref{teo1} to obtain the stationary distribution
for the driving Markov chain. Numerical values are shown in Table \ref{Table3}.
Geddes and Saul (\citeyear{GS}) obtained an approximate stationary distribution
from their simulation, essentially accurate to three decimal places.
One drawback of a simulation in this context is that it does not
clearly show that, when the stationary distribution is expressed as a
matrix, several entries in each column are equal.

We calculate from (\ref{p^circ:d=1}), (\ref{mu^*}) and (\ref{p^*}) that
\[
p^\circ\approx0.0168011,\qquad  \mu^*\approx0.838811,\qquad
p^*\approx0.171001.
\]
Based on their simulation of 1,000,000 coups, Geddes and Saul (\citeyear{GS})
obtained the estimates $0.016638$, $0.838995$ and $0.171451$,
respectively. They did not attempt to estimate the variance parameter.
Using (\ref{var}) and (\ref{covsum}), we find that
\[
\operatorname{Var}(S_0^*)\approx69.860263, \qquad \sum_{m=1}^\infty
\operatorname{Cov}
(S_0^*,S_m^*)\approx-0.951088,
\]
hence
\[
(\sigma^*)^2\approx6.795809.
\]
All displayed numbers are exact except for rounding.

Geddes and Saul (\citeyear{GS}) also proposed a very interesting betting
strategy: Simply play the machine until, and only until, a payout
occurs. Let $E(i,j)$ be the player's expected profit when starting from
cam position $i$ and pointer position $j$. Then
\[
E(i,9)=-1+\mu_i+10q_i, \qquad  i=0,1,\ldots,9,
\]
where of course the 10 is the Futurity award. Furthermore,
\[
E(i,j)=-1+\mu_i+q_iE\bigl(i+1\mbox{ } (\operatorname{mod} 10),j+1\bigr), \qquad  i=0,1,\ldots,9,
\]
for $j=8,7,\ldots,0$ (in that order). These expectations are evaluated
numerically in Table \ref{Table4}. This result is due to Geddes and Saul.

\begin{sidewaystable}
\tablewidth=\textheight
\tablewidth=\textwidth
\caption{Expected player profit when playing until a
payout occurs, as a function of initial cam position (row) and pointer
position (column), rounded to six decimal places; columns 8 and 9 are exact}\label{Table4}
\begin{tabular*}{\textwidth}{@{\extracolsep{\fill}}ld{3.6}d{3.6}d{3.6}ccccccc@{}}
\hline
& \multicolumn{1}{c}{\textbf{0}} & \multicolumn{1}{c}{\textbf{1}} & \multicolumn{1}{c}{\textbf{2}} & \textbf{3} & \textbf{4} & \textbf{5} & \textbf{6} & \textbf{7} & \textbf{8} & \textbf{9} \\
\hline
0 & -1.640567 & -0.210554 & 0.085131 & 0.390591 & 0.706148 &5.122320 & 6.035454 & 6.978775 & 7.953280 & 8.960 \\
1 & -1.056559 & -0.950999 & 0.526288 & 0.831747 & 1.147305 &1.473294 & 6.035454 & 6.978775 & 7.953280 & 8.960 \\
2 & -0.453244 & -0.347685 & -0.238636 & 1.287487 & 1.603045 &1.929034 & 2.265799 & 6.978775 & 7.953280 & 8.960 \\
3 & 0.170015 & 0.275574 & 0.384623 & 0.497277 & 2.073850 & 2.399839 &2.736605 & 3.084503 & 7.953280 & 8.960 \\
4 & 0.813877 & 0.919437 & 1.028486 & 1.141140 & 1.257518 & 2.886209 &3.222975 & 3.570873 & 3.930272 & 8.960 \\
5 & 1.479024 & 1.584584 & 1.693633 & 1.806287 & 1.922665 & 2.042890 &3.725423 & 4.073321 & 4.432720 & 4.804 \\
6 & -0.743671 & 0.686343 & 0.982028 & 1.287487 & 1.603045 & 1.929034 &2.265799 & 6.978775 & 7.953280 & 8.960 \\
7 & -0.130013 & -0.024453 & 1.452834 & 1.758293 & 2.073850 &2.399839 & 2.736605 & 3.084503 & 7.953280 & 8.960 \\
8 & 0.503931 & 0.609491 & 0.718540 & 2.244663 & 2.560220 & 2.886209 &3.222975 & 3.570873 & 3.930272 & 8.960 \\
9 & 1.158832 & 1.264392 & 1.373441 & 1.486095 & 3.062668 & 3.388657 &3.725423 & 4.073321 & 4.432720 & 4.804 \\
\hline
\end{tabular*}
\end{sidewaystable}

We find that, if the pointer position is 3--9, a positive expectation
is assured (regardless of the cam position). In fact, 90 of the 100
expectations are positive. Perhaps more surprising is the fact that
\[
\sum_{(i,j)\in\Sigma} \pi(i,j)E(i,j)\approx0.960501.
\]
In other words, the ``stop after the next payout'' betting system has
positive expectation when played at equilibrium. This observation,
however, is less useful than it may first appear to be. For if the
player has reached approximate equilibrium through extensive play, then
the positive expected profit the system promises will not make up for
the negative expected profit already incurred. And the player should
not expect to find a machine at approximate equilibrium after extensive
play by others. Indeed, a player quits not at a fixed time, such as
after the 10,000th coup, but rather at a random stopping time, such as
after the next win, or after running out of coins. Moreover, as we have
seen, if the pointer position is 3--9, a~player has positive equity and
may not want to relinquish it by walking away. It seems likely that
most players would notice this at least for pointer positions~7, 8~and
9, for in those cases a loss is impossible.

Geddes and Saul (\citeyear{GS}) remarked that ``the machine tends to leave the
player at an unprofitable starting point most of the time after paying
off.'' One way to confirm this is to evaluate the asymptotic
distribution of the Markov chain's state after a payout. Arguing as in
(\ref{p^*}), we get
\begin{eqnarray*}
&&\lim_{n\to\infty}{\sum_{l=1}^n1_{\{R_l^*>0, (X_l,Y_l)=(i,0)\}
}\over\sum_{l=1}^n1_{\{R_l^*>0\}}}\nonumber\\
&& \qquad =\lim_{n\to\infty}{(1/n)\sum_{l=1}^n(1_{\{X_{l-1}=i-1,
R_l>0\}}+1_{\{(X_{l-1},Y_{l-1})=(i-1,9), R_l=0\}})\over(1/n)\sum
_{l=1}^n(1_{\{R_l>0\}}+1_{\{Y_{l-1}=9, R_l=0\}})}\nonumber\\
&& \qquad ={(0.1)p_{i-1}+\pi(i-1,9)q_{i-1}\over(0.8)p_\mathrm
{E}+(0.2)p_\mathrm{O}
+p^\circ}=:\rho(i,0)\qquad  \mbox{a.s.},
\end{eqnarray*}
where $p_{-1}:=p_9$, etc. We find that $\rho(0,0)=\rho(6,0)\approx
0.416991$ and $\rho(i,0)\approx0.020752$ otherwise, and of course
states $(0,0)$ and $(6,0)$ have negative entries in Table \ref
{Table4}. Geddes and Saul obtained approximations from their
simulation. Observe that states $(0,0)$ and $(6,0)$ account for about
0.833982 of the probability, which can be interpreted as the long-term
proportion of payouts that occur when the machine is in mode O. This is
the same as the long-term proportion of Futurity awards that occur when
the machine is in mode O.

Finally, we observe that the previous mode (E or O) is clear at a
glance. This depends on the fact that the machine's payout window
displays not only the three symbols on the payline (from the last coup)
but also the symbols on the line above and the line below the payline.
If the previous mode was O, then exactly four coups are needed to
determine the cam position with certainty; if the previous mode was E,
then at least four and at most seven coups are needed. The player who
is unwilling to play without a positive expectation should play with
pointer position 3 or greater, but also pointer position 2 if the
previous mode was O.

\section{A two-armed slot machine}\label{SectionParrondo}

Motivated by Parrondo's paradox, here we consider a two-armed
generalization of the Futurity slot machine, and we label the arms $A$
and $B$. Excluding the Futurity award, the sequence of payouts from
each arm is assumed nonnegative i.i.d., with arm $A$ (resp., $B$)
having probability~$p_A$ (resp., $p_B$) of a nonzero payout and mean
payout $\mu_A$ (resp., $\mu_B$) in all cam positions. The two arms
are linked only by the Futurity award: After $J$ consecutive losses,
regardless of the order of play of the two arms, the $J$ lost coins are
returned to the player. We assume that $J\ge2$, and we let
$q_A:=1-p_A$ and $q_B:=1-p_B$.

The asymptotic mean payout per coup, including the Futurity award, from
playing arm $A$ only (resp., arm $B$ only) is
\[
\mu_A^*=\mu_A+Jp_A^\circ\qquad  \mbox{where }  p_A^\circ:={p_A
q_A^J\over1-q_A^J}
\]
[resp., $\mu_B^*=\mu_B+Jp_B^\circ$, where $p_B^\circ:=p_B q_B^J/(1-q_B^J)$].
If we play arm $A$ with probability $\gamma$ ($0<\gamma<1$) and arm
$B$ otherwise, a strategy we denote by $C:=\gamma A+(1-\gamma)B$, then
this random mixture has probability
$p_C:= \gamma p_A +(1-\gamma) p_B$ of a nonzero payout and mean payout
$\mu_C:=\gamma\mu_A + (1-\gamma) \mu_B$ in all cam positions,
excluding the Futurity award. Let $q_C:=1-p_C$. Then the asymptotic
mean payout per coup, including the Futurity award, from playing the
random-mixture strategy with parameter $\gamma$ is
\begin{eqnarray*}
\mu^*_C := \mu_C + J p_C^\circ \qquad  \mbox{where }  p_C^\circ
:={p_C q_C^J \over1-q_C^J}.
\end{eqnarray*}

We will say that the \textit{Parrondo effect} is present for the
random-mixture strategy with parameter $\gamma$ if
\[
\mu^*_C<\gamma\mu^*_A + (1-\gamma) \mu^*_B.
\]
In words, the asymptotic mean payout per coup from playing the
random-mixture strategy on the two-armed machine is less than the
asymptotic mean payout per coup from playing the same random-mixture
strategy on two one-armed machines, one of them equivalent to arm $A$
and the other equivalent to arm $B$, each with its own Futurity award.

\begin{theorem}\label{teo4}
If $p_A\ne p_B$, $J\ge2$ and $0<\gamma<1$, then the Parrondo effect
is present for the random-mixture strategy with parameter $\gamma$.
\end{theorem}

\begin{remark*}
As Abbott (\citeyear{A}) remarked, ``In its most general form, Parrondo's
paradox can occur where there is a nonlinear interaction of random
behavior with an asymmetry.'' Here $J\ge2$ ensures the nonlinearity,
while $p_A\ne p_B$ ensures the asymmetry.

In the scenario of Pyke (\citeyear{P}) described in Section \ref{sec1}, $\mu_A^*=\mu
_B^*=1$ (both arms are fair), hence $\mu_C^*<1$ (the random
mixture-strategy is losing for the player, hence winning for the casino).
\end{remark*}

\begin{pf*}{Proof of Theorem \ref{teo4}}
The function $f(x):=(1-x)x^J/(1- x^J)$ is strictly convex on $(0,1)$
for each $J \geq2$ because
\begin{eqnarray*}
f''(x)&=&{Jx^{J-2}[J(1-x)(1+x^J)-(1+x)(1-x^J)]\over(1-x^J)^3}\\
&=&{J(1-x) x^{J-2} \over(1-x^J)^3} \sum_{j=1}^{J-1}
(1-x^j)(1-x^{J-j})\\
&>&0,\qquad  0<x<1.
\end{eqnarray*}
Therefore,
\begin{eqnarray*}
&&\mu^*_C-[\gamma\mu^*_A + (1-\gamma) \mu^*_B]\\
&&\qquad =J\{p_C^\circ-[\gamma p_A^\circ+(1-\gamma)p_B^\circ]\}\\
&& \qquad =J\bigl\{f\bigl(\gamma q_A+(1-\gamma)q_B\bigr)-[\gamma f(q_A)+(1-\gamma
)f(q_B)]\bigr\}\\
&& \qquad <0
\end{eqnarray*}
since $q_A\ne q_B$ and $0<\gamma<1$.
\end{pf*}

In fact, the function $f$ of the proof satisfies
\begin{eqnarray*}
1-x+Jf(x)&=&(1-x){1+(J-1)x^J\over1-x^J}={1+(J-1)x^J\over1+x+\cdots+x^{J-1}}<1
\end{eqnarray*}
for $0<x<1$.
In particular, $p_A+Jp_A^\circ<1$ and $p_B+Jp_B^\circ<1$. If we
assume nonnegative integer payouts, then $\mu_A\ge p_A$ and~$\mu_B\ge
p_B$. It follows that $\mu_A$ and~$\mu_B$ can be chosen in such a way
that $\mu_A^*=\mu_B^*=1$. Actually, fractional payouts per unit bet
are commonplace on modern slot machines. (For example, a machine with
five paylines might return three coins from a five-coin bet.) In such
cases, a loss, for the purpose of the Futurity award, means a zero
payout, not just a payout that is less than the amount bet.

Now we turn to strategies involving nonrandom patterns of the two arms.
Let $D$ denote a (finite) nonrandom pattern of $A$s and $B$s, with at
least one $A$ and at least one $B$, that is repeated ad infinitum. For
example, $D$ could be as simple as $AB$ or~$ABB$, or it could be more
complicated, such as $ABBAB$. Let $r\ge1$ and $s\ge1$ be the numbers
of $A$s and $B$s, respectively, in pattern $D$. Then the asymptotic
mean payout per coup, including the Futurity award, from playing
pattern $D$ repeatedly is given by (\ref{mu^*}) with $I$ equal to the
least common multiple of $r+s$ and $J$. More precisely,
\[
\mu_D^*:={r\mu_A+s\mu_B\over r+s}+Jp_D^\circ,
\]
where $p_D^\circ$ can be inferred from (\ref{p^circ}). The simplest
case is that in which $r+s$ divides $J$ because then (\ref
{p^circ:d=1}) applies and we have
%
\begin{equation}\label{p_D^circ}
p_D^\circ:=\frac{r p_A+s p_B}{r+s}\biggl(\frac{(q_A^r
q_B^s)^{J/(r+s)}}{1-(q_A^r q_B^s)^{J/(r+s)}}\biggr).
\end{equation}
In this case, $p_D^\circ$ (and hence $\mu_D^*$) depends on $D$ only
through $r$ and $s$. For example, $p_{AABBB}^\circ=p_{ABBAB}^\circ$
as long as $J$ is a multiple of 5.

We will say that the \textit{Parrondo effect} is present for the
nonrandom-pattern strategy with pattern $D$ (with $r$ $A$s and $s$
$B$s) if
\[
\mu^*_D<{r\mu^*_A+s\mu^*_B\over r+s}.
\]
The interpretation is analogous to that of the random-mixture strategy.

\begin{theorem}\label{teo5}
If $p_A\ne p_B$, $J\ge2$, and $r,s\ge1$, and if $r+s$ divides $J$,
then the Parrondo effect is present for the nonrandom-pattern strategy
with pattern $D$.
\end{theorem}

\begin{remark*}
While $J$ is a characteristic of the machine, the pattern $D$ (and
hence $r$ and $s$) is chosen by the player, so the assumption that
$r+s$ divides $J$ is too restrictive. We believe that this assumption
can be weakened considerably (see below), but it cannot simply be omitted.
\end{remark*}

\begin{pf*}{Proof of Theorem \ref{teo5}}
The presence of the Parrondo effect is equivalent to
\[
p_D^\circ<{r p_A^\circ+s p_B^\circ\over r+s}.
\]
By the arithmetic mean-geometric mean inequality and $q_A\ne q_B$,
\[
(q_A^r q_B^s)^{1/(r+s)}<{r q_A+s q_B\over r+s}.
\]
Since the function $g(x):=x^J/(1-x^J)$ is increasing on $(0,1)$, we have
\begin{eqnarray*}
p_D^\circ&=&\frac{r p_A+s p_B}{r+s}\biggl(\frac{(q_A^r
q_B^s)^{J/(r+s)}}{1-(q_A^r q_B^s)^{J/(r+s)}}\biggr)\\
&<&{r p_A+s p_B\over r+s}\biggl(\frac{[(r q_A+s q_B)/(r+s)]^J}{1-[(r
q_A+s q_B)/(r+s)]^J}\biggr)\\
&=&p_C^\circ<{r p_A^\circ+s p_B^\circ\over r+s},
\end{eqnarray*}
where $p_C^\circ$ is as in the proof of Theorem \ref{teo4} with $\gamma
:=r/(r+s)$, and the second inequality uses Theorem \ref{teo4}.
\end{pf*}

Various attempts have been made at explaining why Parrondo's paradox
holds in the nonrandom-pattern case; see, for example, Ethier and Lee
(\citeyear{EL}). When the assumptions of Theorem \ref{teo5} are met, we have an
especially simple explanation: the~AM-GM inequality and convexity.

Let us generalize (\ref{p_D^circ}) to arbitrary $D$, $r$, $s$ and
$J$. Although we can minimize the number of terms by taking $I$ to be
the least common multiple of $r+s$ and $J$, we can equally well take
$I$ to be \textit{any} multiple of $r+s$ and $J$, and the simplest
choice is $I:=(r+s)J$ (i.e., $d:=r+s$). Define each of $p_1,p_2,\ldots
,p_{r+s}$ to be $p_A$ or $p_B$ in accordance with the corresponding
term in the pattern $D$. Extend this definition by $p_{i+r+s}=p_i$ for
all $i\in\{1,2,\ldots,r+s\}$, and define $q_i:=1-p_i$ for
$i=1,2,\ldots,2(r+s)$. With this notation, we can write
\[
p_D^\circ={1\over r+s}\sum_{k=1}^{r+s}\Biggl(\sum
_{j=1}^{r+s}p_j\prod_{i=j+1}^{j+kJ-(r+s)\lfloor kJ/(r+s)\rfloor
}q_i\Biggr){(q_A^r q_B^s)^{\lfloor kJ/(r+s)\rfloor}\over1-(q_A^r q_B^s)^J},
\]
where empty products are 1. For example, if $r+s$ divides $J$, then all
products are empty and this reduces algebraically to (\ref{p_D^circ}).
For a less trivial example, consider $D=ABB$. Then, if $J=3K+1$ for a
positive integer $K$,
\begin{eqnarray*}
p_{ABB}^\circ&=&(1/3)[(p_Aq_B+p_Bq_B+p_Bq_A)(q_Aq_B^2)^K\\
&& \hspace*{26pt} {}+(p_Aq_B^2+p_Bq_Bq_A+p_Bq_Aq_B)(q_Aq_B^2)^{2K}\\
&& \hspace*{87pt}\hspace*{26pt} {}+(p_A+2p_B)(q_Aq_B^2)^J]/[1-(q_Aq_B^2)^J],
\end{eqnarray*}
and, if $J=3K+2$ for a nonnegative integer $K$,
\begin{eqnarray*}
p_{ABB}^\circ&=&(1/3)[(p_Aq_B^2+p_Bq_Bq_A+p_Bq_Aq_B)(q_Aq_B^2)^K\\
&&\hspace*{27pt}  {}+(p_Aq_B+p_Bq_B+p_Bq_A)(q_Aq_B^2)^{2K+1}\\
&&\hspace*{73pt} \hspace*{27pt} {}+(p_A+2p_B)(q_Aq_B^2)^J]/[1-(q_Aq_B^2)^J].
\end{eqnarray*}

Despite the impression that may be given by the proof of Theorem \ref{teo5}, it
is not true in general that $p_D^\circ<p_C^\circ$ when $\gamma
:=r/(r+s)$, and it is easy to find counterexamples. It is also not true
in general that $p_D^\circ$ depends on $D$ only through $r$ and $s$.
For example, with $J=6$, $p_{AABB}^\circ>p_{ABAB}^\circ=p_{AB}^\circ
$ if $p_A\ne p_B$. However, extensive numerical computation suggests
the following.

\begin{conjecture*}
Under the assumptions of Theorem \ref{teo5}, the conclusion holds for
patterns of the form $D:=A^rB^s$ if we replace the assumption that
$r+s$ divides $J$ by any one of the following four assumptions:
\begin{longlist}
\item[(a)] $J=2$.

\item[(b)] $\min(r,s)=1$.

\item[(c)] $r+s\le J$.

\item[(d)] $p_A+p_B>1/3$.
\end{longlist}
\end{conjecture*}

We can confirm the sufficiency of condition (b) at least in the
simplest case, $r=s=1$. The case of even $J$ is covered by Theorem \ref{teo5},
so we suppose that $J$ is odd, say $J=2K+1$ for some positive integer
$K$. Then, by algebra,
\begin{eqnarray*}
p_{AB}^\circ-{1\over2}(p_A^\circ+p_B^\circ)
&=&{(p_Aq_B+p_Bq_A)(q_A q_B)^K+(p_A+p_B)(q_A q_B)^J\over2[1-(q_A
q_B)^J]}\\
&& {}-{1\over2}\biggl({p_A q_A^J\over1-q_A^J}+{p_B q_B^J\over
1-q_B^J}\biggr)\\
&=&-{h(q_A,q_B) \over2(1-q_A^J)(1-q_B^J)[1-(q_A q_B)^J]} \\
&<&0,
\end{eqnarray*}
where
\begin{eqnarray*}
h(x,y)&:=&[x^{K+1}-y^{K+1}+(xy)^{K+1}(x^K-y^K)] \\
&& {}\cdot[x^K(1-x)(1-y^{2K+1})-y^K(1-y)(1-x^{2K+1})]\\
&=&(1-x)(1-y)[x^{K+1}-y^{K+1}+(x y)^{K+1}(x^K-y^K)] \\
&& {}\cdot\sum_{k=0}^{K-1}(x^{K-k}-y^{K-k})[(xy)^k-(x y)^{K}]\\
&>&0, \qquad  x,y\in(0,1),  x\ne y.
\end{eqnarray*}

We conclude this section by asking, at what rate can the casino make
money with our two-armed machine, assuming that both arms are fair in
the sense that $\mu_A^*=\mu_B^*=1$? For simplicity, we suppose the
player adopts the random-mixture strategy with $\gamma={1\over2}$.
Then the casino's win rate is
\begin{eqnarray}\label{winrate}
 &&J\biggl[{1\over2}(p_A^\circ+p_B^\circ)-p_C^\circ\biggr]\nonumber\\[-8pt]\\[-8pt]
&&\qquad =J\biggl[{1\over2}\biggl({p_Aq_A^J\over
1-q_A^J}+{p_Bq_B^J\over1-q_B^J}
\biggr)-{[(p_A+p_B)/2][(q_A+q_B)/2]^J\over1-[(q_A+q_B)/2]^J}
\biggr],\nonumber
\end{eqnarray}
which for fixed $J\ge2$ has supremum ${1\over2}[1-J
2^{-J}/(1-2^{-J})]$, achieved as $p_A\to0$ and $p_B\to1$ (and vice
versa). But this case is unrealistic.

Kilby, Fox and Lucas (\citeyear{KFL}), page~137, reported a simulation study of the
effect of hit frequency on player longevity. They considered 10 slot
machines with hit frequencies ranging from 6.7\% to 29.6\% and mean
payouts being roughly equal. So we take $p_A=3/10$ and $p_B=1/15$ as
being the extremes among hit frequencies considered typical in the
industry (for single-payline machines). Notice that condition (d) of
the conjecture is met. We find that (\ref{winrate}) is increasing in
$J$ for $J\le20$ and decreasing in $J$ for $J\ge20$. At $J=20$ its
value is about 0.161553 (i.e., 16.2\%), while at $J=10$ its value is
about 0.100383. Similar calculations can be done for other strategies.
It would seem from the numerical evidence that there is a reasonable
profit potential (for the casino) in a two-armed version of the
Futurity with both arms fair and $J=10$. However, it must be recognized
that there are strategies other than those ordinarily associated with
Parrondo's paradox, so our tentative conclusion about the viability of
this machine on the casino floor is premature.

Consider a strategy for which the choice of arm depends on the Futurity
pointer. Specifically, let $K$ and $L$ be positive integers such that
$K+L=J$, and assume that, if the Futurity pointer shows $j$ consecutive
losses and $0\le j\le K-1$, then arm $A$ is pulled, otherwise arm $B$
is pulled. The driving Markov chain has state space $\Sigma_1:=\{
0,1,\ldots,J-1\}$ and one-step transition matrix $\mathbf{P}_1$ defined by
\begin{eqnarray*}
P_1(i,j)=
\cases{
p_A&\quad if $0\le i\le K-1$ and $j=0$,\vspace*{2pt}\cr
q_A&\quad if $0\le i\le K-1$ and $j=i+1$,\vspace*{2pt}\cr
p_B&\quad if $K\le i\le J-2$ and $j=0$,\vspace*{2pt}\cr
q_B&\quad if $K\le i\le J-2$ and $j=i+1$,\vspace*{2pt}\cr
1&\quad if $i=J-1$ and $j=0$.}
\end{eqnarray*}
This chain is irreducible and aperiodic, and its unique stationary
distribution $\bolds{\pi}_1$ is given by
\begin{eqnarray*}
\pi_1(j)=
\cases{
c^{-1}q_A^j&\quad if $0\le j\le K-1$,\vspace*{2pt}\cr
c^{-1}q_A^K q_B^{j-K}&\quad if $K\le j\le J-1$,}
\end{eqnarray*}
where
\[
c:=1+q_A+\cdots+q_A^{K-1}+q_A^K(1+q_B+\cdots+q_B^{L-1}).
\]

If the mean payouts of arms $A$ and $B$ are
\[
1=\mu_A^*=\mu_A+Jp_A^\circ\quad  \mbox{and}\quad  1=\mu_B^*=\mu
_B+Jp_B^\circ,
\]
then the mean payout at equilibrium under our strategy is
\begin{eqnarray*}
\mu^*&=&\Biggl(\sum_{j=0}^{K-1}\pi_1(j)\Biggr)\mu_A+\Biggl(\sum
_{j=K}^{J-1}\pi_1(j)\Biggr)\mu_B+J\pi_1(J-1)q_B\\
&=&1-\Biggl(\sum_{j=0}^{K-1}\pi_1(j)\Biggr)Jp_A^\circ-\Biggl(\sum
_{j=K}^{J-1}\pi_1(j)\Biggr)Jp_B^\circ+J\pi_1(J-1)q_B,
\end{eqnarray*}
and we find that the Parrondo effect (in favor of the casino) holds if
and only if
\[
\pi_1(J-1)q_B<\Biggl(\sum_{j=0}^{K-1}\pi_1(j)\Biggr)p_A^\circ+
\Biggl(\sum_{j=K}^{J-1}\pi_1(j)\Biggr)p_B^\circ.
\]
Now if we substitute the formula for the stationary distribution, the
constant $c$ is irrelevant, and the condition becomes
\[
q_A^K q_B^L<(1+q_A+\cdots+q_A^{K-1}){p_Aq_A^J\over
1-q_A^J}+q_A^K(1+q_B+\cdots+q_B^{L-1}){p_Bq_B^J\over1-q_B^J}
\]
or
\[
q_A^K q_B^L<{(1-q_A^K)q_A^J\over1-q_A^J}+{q_A^K(1-q_B^L)q_B^J\over1-q_B^J}.
\]
This is equivalent to
\begin{eqnarray*}
q_A^K\biggl({1-q_A^L\over1-q_A^J}-{1-q_B^L\over1-q_B^J}\biggr)<0,
\end{eqnarray*}
which holds if and only if $q_A>q_B$ or equivalently $p_A<p_B$. Here we
are using the fact that the function $f_1(x):= (1-x^L)/(1-x^J)$ is
decreasing on $(0,1)$, which follows from
\begin{eqnarray*}
f_1'(x)&=&{J x^{J-1}(1-x^L)-Lx^{L-1}(1-x^J)\over(1-x^J)^2}\\
&=&-{(1-x)x^{L-1}\over(1-x^J)^2} \sum_{k=1}^K\sum_{l=1}^L
x^{k-1}(1-x^{K-k+l})\\
&<&0,\qquad   0<x<1.
\end{eqnarray*}
So we have the Parrondo effect if $p_A<p_B$. If, however, $p_A>p_B$,
then the Parrondo effect fails and the \textit{player} has the advantage.

Returning to our example in which $p_A=3/10$ and $p_B=1/15$, we suppose
that $J=10$ and consider the above strategy with $K=4$. If $\mu
_A^*=\mu_B^*=1$, then we find that the player's win rate is about
0.145747 (i.e., 14.6\%). We conclude that our machine is not ready for
casino play.

\begin{figure}[b]

\includegraphics{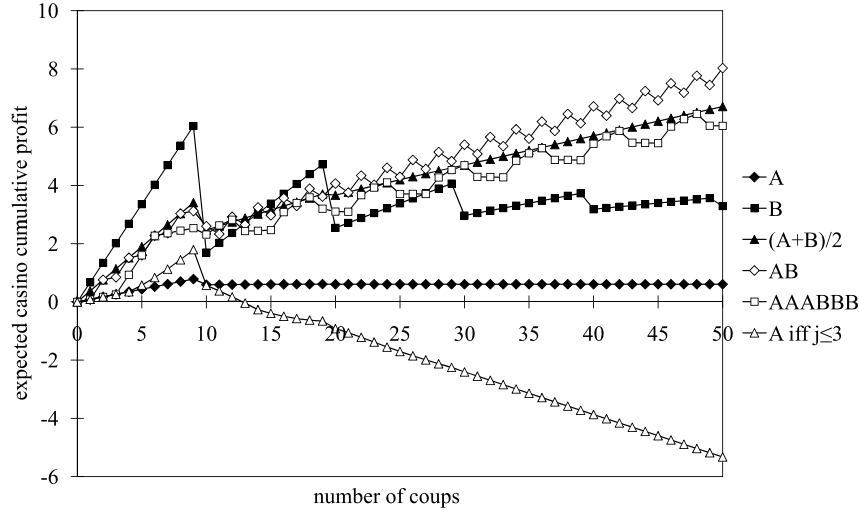}

\caption{Expected casino cumulative profit for various player
strategies. We assume a two-armed slot machine with hit frequencies
$p_A=3/10$ and $p_B=1/15$; Futurity award paid after $J=10$ consecutive
losses, regardless of the order of play of the two arms; initial
pointer position 0; and both arms fair when played exclusively ($\mu
_A^*=\mu_B^*=1$). $j$ is the Futurity pointer position.
Results are by direct calculation (not simulation).}\label{fig1}
\end{figure}
Figure \ref{fig1} compares several strategies in terms of the expected casino
cumulative profit.

The fact that the player can achieve a substantial advantage by using
the information available from the Futurity pointer will not come as a
surprise to those familiar with the original history-dependent Parrondo
games [Parrondo, Harmer and Abbott (\citeyear{PHA})]. Let us recall the
assumptions: In game $A$ the player tosses a \mbox{1$/$2-coin} (heads has
probability 1$/$2), whereas in game $B$, the player tosses a \mbox{9$/$10-coin} if
his last two results are two losses, a \mbox{1$/$4-coin} if his last two results
are a loss and a win in either order, and a \mbox{7$/$10-coin} if his last two
results are two wins. In both games, the player wins one unit with
heads and loses one unit with tails. If the player can use information
about his two most recent results to choose which game to play, the
optimal strategy is clear: Play game $A$ if the last two results differ
and game $B$ otherwise. Most studies of Parrondo's paradox disregard
this strategy and consider only ``blind'' strategies, those that do not
rely on the player's past. In the casino setting, however, one cannot
expect a player to disregard information that may prove to be profitable.

%

\printaddresses

\end{document}